\documentclass[twoside, english, 11pt]{article}

\usepackage[text={6.5in,8.5in},centering]{geometry}
\usepackage{fancyhdr}
\usepackage{natbib}

\usepackage{graphicx} 
\usepackage[textwidth=8em,textsize=small]{todonotes}
\usepackage{amsmath}
\usepackage{natbib}
\usepackage{float}
\usepackage{tikz}
\usepackage{neuralnetwork}
\usepackage{changepage}

\usepackage{enumerate}
\usepackage{url} 
\usepackage{amsfonts,amssymb}
\usepackage{multirow}
\usepackage{booktabs}
\usepackage{tabularx}

\newtheorem{theorem}{Theorem}[section]

\newtheorem{example}{Example}[section]

\newtheorem{corollary}{Corollary}[section]

\newcommand\Author{Zhong \& Wang}
\newcommand\Title{Deep Quantile Regression}

\usepackage{sectsty}
\sectionfont{\centering\large}
\subsectionfont{\bf\normalsize}
\subsubsectionfont{\it\normalsize}

\title{Neural Networks for Partially Linear Quantile Regression}
\date{}
\author {
{Qixian Zhong}\\ Department of Mathematics Science, Tsinghua University\\
Beijing, 100084, China\\ [6mm]
{Jane-Ling Wang}\\ Department of Statistics, University of California, Davis\\
Davis, CA 95616, USA\\ [5mm]
}

\begin{document}
\maketitle 
\begin{abstract}
Deep learning has enjoyed tremendous success in a variety of applications but its application to quantile regressions remains scarce. A major advantage of the deep learning approach is its flexibility to model complex data in a more parsimonious way than nonparametric smoothing methods. However, while deep learning brought breakthroughs in prediction, it often lacks interpretability due to the black-box nature of multilayer structure with millions of parameters, hence it is not well suited for statistical inference. In this paper, we leverage the advantages of deep learning to apply it to quantile regression where the goal to produce interpretable results and perform statistical inference. We achieve this by adopting a semiparametric approach based on the partially linear quantile regression model, where covariates of primary interest for statistical inference are modelled linearly and all other covariates are modelled nonparametrically by means of a deep neural network. In addition to the new methodology, we provide theoretical justification for the proposed model by establishing the root-$n$ consistency and asymptotically normality of the parametric coefficient estimator and the minimax optimal convergence rate of the neural nonparametric function estimator. Across several simulated and real data examples, our proposed model empirically produces superior estimates and more accurate predictions than various alternative approaches.

\noindent%
\vfill
\end{abstract}

\clearpage
\section{Introduction}
\label{sec: introduction}
With advances in computational power and the availability of large data, deep learning has emerged as a  powerful data analysis tool  in a wide variety of applications, such as computer vision \citep{krizhevsky2012imagenet, russakovsky2015imagenet}, speech recognition \citep{hinton2012deep}, and natural language processing \citep{collobert2011natural}. Deep learning estimates maps from data using neural networks which compose of multiple (parameterized) nonlinear transformations. These inferred transformations are jointly optimized \emph{end-to-end} in order to produce the optimal overall map (rather than independently estimating each transformation in a separate stage).

Roughly speaking, a neural network, which consists of several layers and neurons between the input and output layers, is a composite function (see formula \eqref{def: DNN1}) with a recursive concatenation of an affine linear function and a simple nonlinear map.  
The success of neural networks is attributed to their powerful capacity to represent unknown functions. 
For example, \cite{cybenko1989approximation} and \cite{hornik1989multilayer} 
showed that any continuous functions can be approximated by shallow neural networks to any degree of accuracy. \citet{telgarsky2016benefits} and \citet{yarotsky2017error} further showed that deep neural networks enjoy a better representational power than their shallow counterparts. 

Despite their superior empirical performance, deep learning models, mostly a black box, often lack {intepretability and theoretical support.  Different approaches have emerged in recent    works  to examine various aspects of interpretable deep learning models. 
For instance,   saliency-based  \citep{zeiler2014visualizing,simonyan2014deep,selvaraju2017grad}  and concept-based \citep{kim2018interpretability,yeh2020completeness} methods aim at providing post hoc explanations for a certain type of neural networks. Another approach by \cite{chen2018looks} and \cite{li2018deep} focus on designing specific neural network structures for case-based reasoning. Neural networks  have also been adapted to study the causal effects between variables  \citep{luo2020causal,farrell2021deep,shi2019adapting}. For additional   works on intepretable deep learning models, we refer readers to the recent review papers \citep{chakraborty2017interpretability,murdoch2019interpretable,rudin2019stop} and reference therein.

Unlike the above approaches, 
this paper adopts the statistical model-based approaches for interpretability by constructing  neural networks for a partially linear quantile regression (PLQR) problem. Specifically, we model the the covariates of interest with a linear predictor for interpretability and statistical inference and  model the nonparametric component with neural networks. The proposed deep learning method for PLQR is abbreviated as DPLQR. As a  semiparametric  approach, 
DPLQR not only offers interpretibility for the parametric component but also allows model flexibility for the nonparametric component. Importantly,  it avoids the curse of dimensionality of nonparametric smoothing methods through the strength of neural networks to detect the  structure, often low-dimensional, of the data.  
We further provide mathematical support for the DPLQR, which not only quantifies the uncertaity of the inference but somewhat reveals the success of the deep learning.} 

Since the seminal work of \cite{koenker1978regression}, 
quantile regression has been extensively investigated, including linear quantile regression \citep{koenker1978regression,portnoy1991asymptotic}, 
nonparametric quantile regression \citep{samanta1989non,jones1990mean,chaudhuri1991nonparametric,he1994convergence}  
and semiparametric quantile regression \citep{he1996bivariate,lee2003efficient,wu2010single,cai2012semiparametric}.  
For a comprehensive introduction of quantile regression, we refer to the monographs by \cite{koenker2005quantile} and \cite{koenker2017handbook}.  
Compared to the least squares regression approach that focuses on the conditional mean of the response, 
quantile regression offers a more expansive view of the effect of covariates on a response. 
Moreover,  quantile regression is more robust against   outliers  when the distribution of the response is heavy-tailed or skewed.

While linear and nonparametric quantile regression have been well developed, theory and methodology for partially linear quantile regression models are lagging and existing work is mainly focused on the partially linear additive quantile regression  \citep{lian2012semiparametric, hoshino2014quantile, sherwood2016partially}. 
This  approach  incorporates a linear regression for  some covariates  and an additive model with  smooth but unknown  regression functions for the remaining covariates. The additive structure alleviates the  curse of dimensionality but it is not amenable to model  interactions among covariates. 
Meanwhile, exisitng  fully nonparametric approaches suffer from a severe curse of dimensionality, so they are only effective for very low dimensional covariates.  To fill these gaps, we consider DPLQR, which models some covariate effects with a linear model but the rest with an unknown  multivariate continuous function.  This model is effective in  interpreting the  effects of primary  covariates, 
 such as the effect of a treatment. It  also enjoys the flexibility of a fully nonparametric function but is more resilient to the curse of dimensionality. Our theoretical results are in line  with recent studies \citep{petersen2018optimal,bauer2019deep,schmidt2020nonparametric} which show that  deep learning has the ability to learn the unknown underlying {low dimensional structure of the data embedded in high dimension space.} 
 This is a major advantage over the traditional smoothing approaches that were designed  to estimate the  covariate effects nonparametrically.

Applications of deep learning to quantile regression have emerged in recent years, such as in climate prediction 
\citep{hatalis2017smooth} 
and 
electricity and power system 
\citep{gan2018embedding}. However, theoretical understanding of quantile regression with neural networks remains scarce and limited to nonparametric quantile regression. 
\cite{romano2019conformalized} employed conformal methods to construct prediction intervals for the response but did not  address  estimation  of the conditional quantile function. 
\cite{jantre2020quantile} developed  consistency results for nonparametric quantile function estimator  with a single-hidden-layer neural network. {However, the implementation of their procedure requires exponential time to compute as compared to the polynomial time  for deep neural networks  \citep{rolnick2017power}.} 
As we were writing up the results of our research findings,  we became aware of a related work that was independently developed by \cite{padilla2020quantile}. Although this work also explored the convergence rate of the conditional quantile function estimator, 
it is substantially different from ours. First, it focuses on a black-box nonparametric approach to estimate  the quantile function, 
while we are interested in both estimation and  interpretability as well as statistical inference for the model. 
Second, the theoretical analysis of their work only holds for continuous covariates while our theory covers both continuous and discrete covariates with asymptotic normality established for the estimates of the linear component. 

To summarize, the major contributions of this paper are four-fold.  
\begin{enumerate}
    \item { We introduce  DPLQR aiming to shed new light on an interpretable deep learning model to overcome the drawback of a black-box deep learning approach. {Although there are a number of attempts to address it,} most of them fail to provide uncertainty quantification.  
    In contract, we develop confidence intervals for the effects of linear  covariates, which  are of interest to  practitioners.}   
     Our approach can thus be viewed as a bridge between machine  learning and statistical inference. 
    \item We provide  theoretical justification for deep learning research by showing the minimax optimal convergence rates (up to a poly-logarithmic factor) of the nonlinear component of the DPLQR. We further establish  asymptotic normality of the  regression coefficient estimator for  both homoscedastic and heteroscedastic random errors. 
    \item The proposed DPLQR model  is flexible and includes  a large number of previously-studied quantile regression models. Specifically,  DPLQR  reduces to linear quantile regression when the nonparametric component is absent and  it reduces to  nonparametric quantile regression in the absence of linear predictors. The DPLQR model also includes the partially linear additive quantile regression model. 
    
\item    Our methodology is able to identify the underlying intrinsic dimension of the data, which circumvents the  curse-of-dimensionality incurred by a nonparametric smoothing approach.   For example, when the true model corresponds to a  partially linear additive quantile regression, 
    the resulting neural network estimators have one-dimensional nonparametric rates of convergence (up to a poly-logarithmic factor). 
\end{enumerate}

The rest of the paper proceeds as follows. In section \ref{sec: preliminaries}, we briefly introduce the fundamental concept of neural networks and quantile regression. Asymptotic properties of the estimators are presented in Section \ref{sec: results}. The  implementation of the proposed approach is discussed in Section \ref{sec: computational details} 
along with the calculation of the asymptotic covariance matrix for the vector parameter. Section \ref{sec: simulations} and Section \ref{sec: application} provide simulation studies and data applications comparing the proposed method with linear quantile regression and partially linear additive quantile regression. 
Section \ref{sec: funture works} discusses some potential extensions { and Section \ref{sec: proofofthm} provides proofs of the theorems}. 

\section{Preliminaries}
\label{sec: preliminaries}

\subsection{Neural network}
\label{subsec: nueral network}
\begin{figure}
\centering
\includegraphics[width=9.0cm,height=5cm]{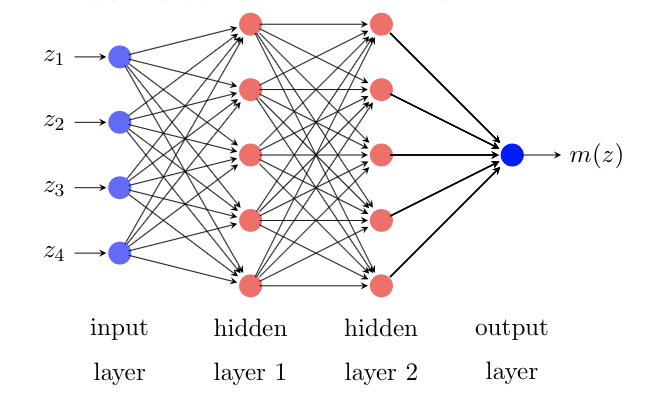}
\caption{\label{fig: network} A 3-layer neural network with four input variables and one output.}
\end{figure}

We first briefly present relevant  background on deep neural networks. For some integer $L\ge 2$, let  $\boldsymbol{q}=(q_0,q_1,\ldots,q_L)^\top\in\mathbb{N}^{L+1}$.  An $L$-layer neural network with input dimension $q_0$ and output dimension $q_L$ is a function $m:\mathbb{R}^{q_0}\rightarrow \mathbb{R}^{q_{L}}$ that satisfies the following recursive relation: 
\begin{equation}\label{def: DNN1}
    \begin{aligned}
    & m(z)=\tilde{W}_L m_{L-1}(z)+\tilde{b}_L,\\
    &m_{L-1}(z) = \sigma(\tilde{W}_{L-1}m_{L-2}(z)+\tilde{b}_{L-1}),\\
    &\ldots,\\
     &m_1(z)=\sigma(\tilde{W}_1m_0(z)+\tilde{b}_1),\\&m_0(z) = z,
    \end{aligned}
\end{equation}
where $\tilde{W}_k$ and $\tilde{b}_k$ are $q_{k-1}\times q_k$ matrix and $q_k$-dimensional column vector, respectively, and $\sigma$ is a prior deterministic function which operates component-wise on vectors, i.e., 
$\sigma((v_1,\ldots,v_m)^\top)=(\sigma(v_1),\ldots,\sigma(v_m))^\top.$ We call $L$ the \textit{depth} of the neural network, $m_k$ for $1\le k\le L-1$ the $k$-th \textit{hidden layer} and $\sigma:\mathbb{R}\rightarrow\mathbb{R}$ the \textit{activation function}. Two layers ($L=2$) one is often called a \textit{shallow neural network}. 
At the $k$-th hidden layer, there are $q_k$ neurons, or nodes, and  $q_k$ is called the  \textit{width} of the neural network. The activation function $\sigma$ links adjacent layers and  is often set to be a simple nonlinear function. In this paper, we consider the \textit{rectified linear unit} (ReLU) activation function $\sigma(z)=\max(z,0)$ since it is computationally efficient and often achieves best  practical performance in practice \citep{krizhevsky2012imagenet}. The matrices $\tilde{W}_k$ and vectors $\tilde{b}_k$  are often referred to as the ``weight''  and ``bias'' respectively in the machine learning  literature, but  we avoid using these terms here to prevent  confusion.  We denote $W_k=(\tilde{W}_k, \tilde{b}_k)\in\mathbb{R}^{q_k\times (q_{k-1}+1)}.$ Then the neural network in  \eqref{def: DNN1} can be succinctly expressed as 
\begin{equation}\label{def: DNN2}
    m(z)=W_{L}\tilde{\sigma}\circ\cdots\circ W_2\tilde{\sigma}(W_1 \tilde{z}), 
\end{equation}
where $\tilde{\sigma}(v)=(\sigma(v)^\top,1)^\top$ and $\tilde{z}=(z^\top,1)^\top$. Figure \ref{fig: network} illustrates a three layers neural network with $\boldsymbol{q}=(4,5,5,1)^\top$
  
Note that the total number of parameters in \eqref{def: DNN2} is $\sum_{k=1}^{L}q_{k}(q_{k-1}+1)$, which can be very large and may lead to overfitting. \cite{han2015learning,bauer2019deep} and \cite{schmidt2020nonparametric} mitigated against this  by deactivating  some of the links of neurons between the adjacent hidden layers.
Following this strategy, for $s\in \mathbb{N}$, $L\ge 2$, $A> 0$ and $\boldsymbol{q}=(q_0,q_1,\ldots,q_L)^\top$, we consider a sparsely connected neural network class
\begin{equation}\label{def: sparse NN}
\begin{aligned}
    \mathcal{M}(s,L,\boldsymbol{q},A)=\Big\{ m(z)=&W_{L}\tilde{\sigma}\circ\cdots\circ W_2\tilde{\sigma}(W_1 \tilde{z})~|~ W_k\in\mathbb{R}^{q_k\times(q_{k-1}+1)},~  \|W_k\|_{\infty}\le 1~\text{for}~ \\
    & k =1,\ldots,L,~ \sum_{k=1}^{L}\|W_k\|_{0}\le s ~\text{and}~ \|m\|_{\infty}\le A    \Big\},
\end{aligned}
\end{equation}
where $\|\cdot\|_{\infty}$ is the sup-norm of a matrix or function and $\|\cdot\|_{0}$ is the number of non-zero elements of a matrix.  

\subsection{Partially linear quantile regression model and estimation}
Consider a univariate random variable $Y$ and a multivariate random variable $U=(X,Z)\in \mathbb{R}^{p}\times\mathbb{R}^q$, 
of which $X$ can include treatment variables and continuous covariates of interest. Let $F_{Y|U}(\cdot|u)$ be the conditional distribution function of $Y$ given $U=u.$ For some $0<\tau<1$, the $\tau$-th conditional quantile of $Y$ given $U=u$ is defined as
\begin{equation*}
    \xi_{\tau}(u)=\mathop{\inf}_{y\in\mathbb{R}}\{y~|~F_{Y|U}(y|u)\ge \tau\}.
\end{equation*}

In this paper, we assume $ \xi_{\tau}(X,Z)=X^\top\theta_\tau + m_\tau(Z)$, which leads to the following partially linear quantile regression model: 
\begin{equation}\label{def: quantile model}
    Y=X^\top\theta_\tau + m_\tau(Z)+\epsilon, P(\epsilon\le0|U)=\tau,
\end{equation}
where $\theta_\tau\in\mathbb{R}^p$ is an unspecified parameter, $m_\tau:\mathbb{R}^{q}\rightarrow\mathbb{R}$ is an unknown function and the  error $\epsilon$ can  be heteroscedastic by allowing it to vary with $u=(x,z)$. 

Let $\{(X_i,Z_i,Y_i):i=1,\ldots,n\}$ denote independent and identically distributed realizations of  $(X,Z,Y)$. For simplicity, we use the notation $\mathcal{M}$ to denote the neural network class $\mathcal{M}(s,L,\boldsymbol{q},\infty)$ in (\ref{def: sparse NN}) with $q_0=q$ and $q_L=1$. 
To estimate the vector $\theta_\tau$ and the function $m_\tau$, we minimize the loss function:
\begin{equation}\label{objective function}
    (\hat{\theta}_\tau,\hat{m}_\tau)=\mathop{\arg\min}_{(\theta,m)\in \mathbb{R}^{p}\times \mathcal{M}_{}} \frac{1}{n}\sum_{i=1}^n\rho_\tau(Y_i-X_i^{\top}\theta-m(Z_i)),
\end{equation}
where $\rho_\tau(t)=t\{\tau-1(t<0)\}$ is called the check loss. This loss function becomes the absolute value $L^1$-loss when $\tau=0.5$ which leads to the median estimators. For brevity, we suppress the subscript $\tau$  and write $(\theta_0,m_0)=(\theta_\tau,m_\tau)$ and $(\hat{\theta},\hat{m})=(\hat{\theta}_\tau,\hat{m}_\tau)$. 

\section{Theory}
\label{sec: results}
In this section, we establish the theoretical properties of the estimators $\hat{\theta}$ and $\hat{m}$. 
We first introduce a class of smooth functions in which $m_0$ resides. 

Let $\gamma$ and $B$ be two positive constants and $\lfloor \gamma \rfloor$ denote the largest integer strictly less than $\gamma$. We call a function $h: \mathbb{T}\subset \mathbb{R}^q\rightarrow\mathbb{R}$ a  $(\gamma,B)$-H{\"o}lder smooth function if it satisfies 
\begin{equation*}
\sup_{z\in \mathbb{T}}    \Big|\frac{\partial^{|\boldsymbol{\alpha}|}h}{\partial^{\alpha_1}z_1\ldots\partial^{\alpha_q}z_q}(z)\Big|\le B,~\text{for all}~ \boldsymbol{\alpha}=(\alpha_1,\ldots,\alpha_q)^\top\in\mathbb{N}^q ~\text{and}~|\boldsymbol{\alpha}|=\sum_{i=1}^{q}\alpha_i\le \lfloor \gamma \rfloor,
\end{equation*}
and
\begin{equation*}
   \mathop{\sup}_{z,z^*\in\mathbb{T}} \Big|\frac{\partial^{|\boldsymbol{\alpha}|}h}{\partial^{\alpha_1}z_1\ldots\partial^{\alpha_q}z_q}(z )-\frac{\partial^{|\boldsymbol{\alpha}|}h}{\partial^{\alpha_1}z_1\ldots\partial^{\alpha_q}z_q}(z^*)
    \Big|\le B\|z-z^*\|_{2}^{\gamma-\lfloor \gamma \rfloor},~\text{for all}~|\boldsymbol{\alpha}|=\lfloor \gamma \rfloor.
\end{equation*}

Denote the class of all such $(\gamma,B)$-H{\"o}lder smooth functions as $\mathcal{H}_{q}^{\gamma}(\mathbb{T},B)$. Let 
$J\in \mathbb{N}$, $\boldsymbol{\gamma}=(\gamma_1,\ldots,\gamma_{J})^\top\in\mathbb{R}_{+}^{J}$,  $\boldsymbol{d}=(q,d_1,\ldots,d_{J})^\top\in\mathbb{N}^{J+1}$ and $\boldsymbol{\bar{d}}=(\bar{d}_1,\ldots,\bar{d}_{J})^\top\in\mathbb{N}^{J}$ with $\bar{d}_1\le q$ and $\bar{d}_k\le d_{k-1}, k=2,\ldots,J$. We further define a composite function class:
\begin{equation}\label{def: calH}
\begin{aligned}
        \mathcal{H}(J,\boldsymbol{\gamma},{\boldsymbol{d}},\boldsymbol{\bar{d}},B)=\Big\{h=&h_{J}\circ\ldots\circ h_1:\mathbb{T}\rightarrow\mathbb{R}~|~ \   h_k=(h_{k1},\ldots,h_{k d_k})^\top~\text{and}~\\& h_{kj}\in \mathcal{H}_{\bar{d}_k}^{\gamma_k}([a_k,b_k]^{\bar{d}_k},B) ~\text{for some}~|a_k|,|b_k|\le B \Big\}.
\end{aligned}
\end{equation}
Note that this class of functions, first proposed by \cite{schmidt2020nonparametric}, contains two kinds of dimension ${\boldsymbol{d}}$ and $\boldsymbol{\bar{d}}.$ We call   $\boldsymbol{\bar{d}}$ the \textit{intrinsic dimension} of the function $h$ in $\mathcal{H}(J,\boldsymbol{\gamma},{\boldsymbol{d}},\boldsymbol{\bar{d}},B)$. 

For an illustration,  consider the function
\begin{equation}\label{def: nn7}
\begin{aligned}
            h(z)=h_{31}( h_{21}(h_{11}(z_1,z_2),h_{12}(z_3,z_4)),h_{22}(h_{13}(z_5,z_6),h_{14}(z_7))),
\end{aligned}
\end{equation}
where all $h_{ij}$ are $(\gamma,1)$-H{\"o}lder smooth. It is clear that $h\in\mathcal{H}(J,\boldsymbol{\gamma},{\boldsymbol{d}},\boldsymbol{\bar{d}},B)$ with $J=3,\boldsymbol{\gamma}=(\gamma,\gamma,\gamma)^\top,{\boldsymbol{d}}=(7,4,2,1)^\top,\boldsymbol{\bar{d}}=(2,2,2)^\top$ and $B=1.$

With  different choices of $J,\boldsymbol{\gamma},$ ${\boldsymbol{d}}$ and $\boldsymbol{\bar{d}}$,  $\mathcal{H}(J,\boldsymbol{\gamma},{\boldsymbol{d}},\boldsymbol{\bar{d}},B)$ includes a large number of function classes that have been considered in the statistical and economic literature. Below we provide two examples to illustrate the ubiquity of such function classes.  We say a function $h$ is $(\infty,B)$-H{\"o}lder smooth if it is $(\gamma,B)$-H{\"o}lder smooth for all $\gamma>0$.
\begin{example} [Generalize additive functions]
A function $h:\mathbb{R}^q\rightarrow\mathbb{R}$ is additive if it can be represented a  sum of univariate functions of each components \citep{stone1985additive}, i.e., for $z=(z_1,\ldots,z_q)^\top$,
\begin{equation}\label{mod: additive}
    h(z)= h_1(z_1)+\ldots+h_q(z_q),
\end{equation}
where $h_k, k=1,\ldots,q$ are univariate $(\gamma,B)$-H{\"o}lder smooth  functions. Here $J=2, \boldsymbol{\gamma}=(\gamma,\infty)^\top,$ $\boldsymbol{d}=(q,q,1)^\top,$ $\boldsymbol{\bar{d}}=(1,q)^\top,$ $h_{1k}(z)=h_k(z_k),~k=1,\ldots,q$, and $h_{21}(y)=y_1+\ldots+y_q,$ where $y=(y_1,\ldots,y_q)^\top.$ Furthermore, \cite{horowitz2001nonparametric} added an unknown link function $g$ and proposed the generalized additive function:
\begin{equation*}
    h(z)=g(h_{1}(z_1)+\ldots+h_{q}(z_q)), 
\end{equation*}
where $g$ and $h_k,k=1,\ldots,q$ are univariate $(\gamma,B)$-H{\"o}lder smooth  functions.
In this case, the function $h$ has a hierarchical structure with  $J=3, \boldsymbol{\gamma}=(\gamma,\infty,\gamma)^\top,$ $\boldsymbol{d}=(q,q,1,1)^\top,$ and $\boldsymbol{\bar{d}}=(1,q,1)^\top,$ $h_{1k}(z)=h_k(z_k),~k=1,\ldots,q$,   $h_{21}(y)=y_1+\ldots+y_q,$ for $y=(y_1,\ldots,y_q)^\top,$ and $h_{31}=g.$  


\end{example} 

\begin{example} [Single/multiple index functions]
A single index function, first introduced by \cite{ichimura1991semiparametric} and later extended to multiple indices  by  \cite{hristache2001structure}, is given by :
\begin{equation}
    h(z)=h_1(z^\top\alpha_1,\ldots,z^\top\alpha_K),
\end{equation}
where  $\alpha_{k},k=1,\ldots,K$ are unknown parameters and $z^\top\alpha_j$ are the index functions.
It is easy to see that $h_{1k}(z)=z^\top\alpha_k,k=1,\ldots,K$ and $h_{21}(y)=h_1(y)$.  Thus, if $h_1$ is $(\gamma,B)$-H{\"o}lder smooth,  $\boldsymbol{\gamma}=(\infty,\gamma)^\top,$ $\boldsymbol{d}=(q,K,1)^\top$ and $\boldsymbol{\bar{d}}=(\bar{d}_1,K)^\top$ with $\bar{d}_1=\max_{k} \{\|\alpha_k\|_0\}.$  

\end{example}

 For some $J\in \mathbb{N}$, $\boldsymbol{\gamma}=(\gamma_1,\ldots,\gamma_{J})\in\mathbb{R}_{+}^{J}$,  $\boldsymbol{d}=(q,d_1,\ldots,d_{J})^\top\in\mathbb{N}^{J+1}$ and $\boldsymbol{\bar{d}}=(\bar{d}_1,\ldots,\bar{d}_{J})^\top\in\mathbb{N}^{J}$ with $\bar{d}_1\le q$ and $\bar{d}_k\le d_{k-1}, k=2,\ldots,J ,$  we define the \textit{effective smoothness}  $\bar{\gamma}_{k}= \gamma_{k}\prod_{i=k+1}^{J}(\gamma_i\wedge 1)$  of a function $h$ in $\mathcal{H}(J,\boldsymbol{\gamma},{\boldsymbol{d}},\boldsymbol{\bar{d}},B),$  and write
\begin{equation*}
    \bar{k}=\mathop{\arg\min}_{k\in\{1,\ldots,J\}} \frac{\bar{\gamma}_{k}}{2\bar{\gamma}_{k}+\bar{d}_k}~\text{and}~r_n=n^{-\frac{\bar{\gamma}_{\bar{k}}}{2\bar{\gamma}_{\bar{k}}+\bar{d}_{\bar{k}}}}.
\end{equation*}
For the covariate $X=(X_1,\ldots,X_p)^\top$, we define 
\begin{equation}\label{def: projection}
    \varphi_{k}^{*}=\mathop{\arg\min}_{\varphi\in L^{2}(P_{Z})}\mathbb{E} [f(0|U)\{X_k-\varphi(Z)\}^2], k=1,\ldots,p,
\end{equation}
where $L^{2}(P_{Z})=\{\varphi~|~ \mathbb{E}\varphi^2(Z) < \infty \}.$ And denote  
$\boldsymbol{\varphi}^*(Z)=(\varphi_1^*(Z),\ldots,\varphi_p^*(Z))^\top,$  $\Sigma_1=\mathbb{E}[\tau(1-\tau)\{X-\boldsymbol{\varphi}^*(Z)\}\{X-\boldsymbol{\varphi}^*(Z)\}^\top]$ and  $\Sigma_2=\mathbb{E}[f(0|U)\{X-\boldsymbol{\varphi}^*(Z)\}\{X-\boldsymbol{\varphi}^*(Z)\}^\top].$ It is easy to show that $\boldsymbol{\varphi}^{*}=\mathbb{E}(X|Z)$,  if the conditional error density $f(\cdot|U)$ is independent of $U$ at zero, 
see also \cite{lian2012semiparametric} and \cite{hoshino2014quantile} for partially linear additive regression. 

Next, we state the assumptions for the deep partially linear quantile regression model. 
\begin{enumerate}[{({A}1)}]
\item The true vector parameter $\theta_0$ belongs to a compact subset   $\Theta\subset\mathbb{R}^p$ and the true nonparametric function $m_0$ belongs to $\mathcal{H}=\mathcal{H}(J,\boldsymbol{\gamma},{\boldsymbol{d}},\boldsymbol{\bar{d}},B).$  \label{assump: param space} 
\item The covariates   $(X,Z)$ take values in a compact subset of $\mathcal{R}^{p+q}$ that, without loss of generality, will be assumed to be $[0,1]^{p+q}.$ In addition, the probability density function (PDF) of $Z$ is bounded away from zero and from infinity. \label{assump: covariates}
\item The conditional PDF $f(\cdot|u)$   of the random error $\epsilon$ given the covariate $U=u$, has continuous derivative $f^{'}(\cdot|u),$ and there exist positive constants $b_0$ and $c_0$ 
such that $1/c_0<f(t|u) < c_0$ and $|f^{'}(t|u)|<c_0$ for all $|t|\le b_0, u\in[0,1]^{p+q}.$\label{assump: error pdf}
\item $L =O(\log n),$ $s=O(n r_n^2 \log n)$ and   $n r_n^2\lesssim \text{min}_{k=1,\ldots,L}\{q_k\}\le\text{max}_{k=1,\ldots,L}\{q_k\} \lesssim n.$ \label{assump: neural struc}
\item The matrices $\Sigma_1$ and $\Sigma_2$ are both positive definite. \label{assump: matrix}
\item $\bar{\gamma}_{\bar{k}}>\bar{d}_{\bar{k}}/2$ and { $\max_{k=1,\ldots,p}(\mathbb{E}|X_k|^{4})<\infty.$}  \label{assump: AsympNormal} 
\end{enumerate}

The boundedness of both the parameters and covariate spaces in assumptions (A\ref{assump: param space}) and (A\ref{assump: covariates}) are standard for semiparametric/nonparametric regression. In  (A\ref{assump: error pdf}) we assume that the PDF of the error and its derivative are bounded to guarantee that the true parameter $(\theta_0,m_0)$ is a well-separated point of the minimum of the expected check loss function. For (A\ref{assump: neural struc}), we assume that the size of neural networks $\mathcal{M}$ used in \eqref{objective function} grows with the sample size $n$ at a certain rate to balance the approximation and estimation errors of the estimators.  Assumptions (A\ref{assump: matrix}) and (A\ref{assump: AsympNormal}) are common conditions for asymptotic normality of the vector estimator $\hat{\theta}$ in semiparametric regression \citep{horowitz2009semiparametric}, where (A\ref{assump: matrix}) is used to develop the asymptotic variance while (A\ref{assump: AsympNormal}) guarantees  $\sqrt{n}$-consistency. 

We are now ready to state the convergence rate of the estimators.
\begin{theorem}\label{Thm: convergency rate}
Under Assumptions (A\ref{assump: param space})-(A\ref{assump: matrix}), 
we have
\begin{equation*}
 	    \lim_{C\rightarrow\infty}\lim_{n\rightarrow\infty}\sup_{m_0\in\mathcal{H}}\mathbb{P}(\|\hat{m}-m_0\|_{L^2([0,1]^{q})}\ge C r_n \log^2 n) = 0.
 	\end{equation*}
\end{theorem}


From the proof of Theorem \ref{Thm: convergency rate} one can see that the convergence rate is the result of a trade-off between  estimation error and approximation error. Here the approximation error is defined as the distance between the true parameter $m_0$ and the neural network set $\mathcal{M},$ i.e., $\mathop{\min}_{m\in\mathcal{M}}\|m-m_0\|_{L^2([0,1]^q)}$. 
It is known that a more complex neural network structure is more flexible and thus leads to a smaller approximation error  \citep{anthony2009neural,yarotsky2017error,bauer2019deep,schmidt2020nonparametric}. However, too many parameters will lead to  high variance. Hence, there is an implicit  ``bias-variance'' trade-off that is reflected in  the growth of  neural networks. 

Note that the convergence rate of the estimator $\hat{m}$ is determined by both the effective smoothness and the intrinsic dimension of the true function $m_0$, rather than the dimension $q$ of the covariate $Z$.  For example, if $m_0$ has the composite structure in \eqref{def: nn7},  the convergence rate for the proposed method is $n^{-\gamma/(2\gamma+2)}\log^2 n$. In contrast,  the convergence rate for a nonparametric method, such as  kernel  or spline smoothing   is of the order $n^{-\gamma/(2\gamma+7)}$. 
This shows that our method is able to detect the low dimensional structure of the data and circumvents  the curse of dimensionality.

In particular, when $m_0$ reduces to additive or single index function, 
 the resulting estimators have one-dimensional nonparametric rates of convergence (up to a poly-logarithmic factor). 
 This is similar to  results of \cite{stone1985additive} and \cite{ichimura1991semiparametric} for nonparametric regression. 

The next theorem establishes the minimax lower bound for estimating $m_0$, which implies that the resulting estimator $\hat{m}$ in Theorem \ref{Thm: convergency rate} is rate-optimal.

\begin{theorem}\label{Thm: minimax lower bound}
Let $\mathcal{F}$ be the class of probability density functions that satisfy Assumption (A\ref{assump: error pdf}). 
Then we have 
\begin{equation*}
    \mathop{\lim}_{C\rightarrow\infty}\mathop{\lim}_{n\rightarrow\infty}\mathop{\inf}_{\hat{m}}\mathop{\sup}_{(\theta_0,m_0,f)\in\mathbb{R}^{q}\times\mathcal{H}\times\mathcal{F}}\mathbb{P}_{(\theta_0,m_0,f)}\big(\|\hat{m}-m_0\|_{L^2([0,1]^q)}\ge C r_n\big)=1,
\end{equation*}
where the infimum is taken over all possible predictors $\hat{m}$ based on the observed data.
\end{theorem}

Below we show that the estimator $\hat{\theta}$ for the vector parameter is asymptotically normal at the $\sqrt{n}$ rate.
\begin{theorem}\label{Thm: asymptotic normality}
Under Assumptions (A\ref{assump: param space})-(A\ref{assump: AsympNormal}), 
we have
\begin{equation*}
    \sqrt{n}(\hat{\theta}- \theta_0)\rightarrow N(0,\Sigma_2^{-1}\Sigma_1\Sigma_2^{-1}).
\end{equation*}
\end{theorem}

When $f(0|U)$ is a constant function, the solution of \eqref{def: projection} would be $\boldsymbol{\varphi}^*(Z)=\mathbb{E}(X|Z),$ which  leads to $\Sigma_1=\tau(1-\tau)\text{Var}\{X-\mathbb{E}(X|Z)\}$,  $\Sigma_2=f(0)\text{Var}\{X-\mathbb{E}(X|Z)\}$ and more generally, the following corollary. 

\begin{corollary}\label{cor: indp theta asymp}
Under the same assumptions of Theorem \ref{Thm: asymptotic normality} and when 
$f(0|U)$ is a constant function, 
 we have
\begin{equation*}
    \sqrt{n}(\hat{\theta}- \theta_0)\rightarrow N(0,\Sigma),
\end{equation*}
where $\Sigma={\tau(1-\tau)}[\text{Var}\{X-\mathbb{E}(X|Z)\}]^{-1}/{f^2(0)}.$
\end{corollary}

For partially linear quantile regression with homoscedastic error, the random error $\epsilon$ is independent of the covariate $U$, which implies that $f(0|U=u)=f(0),$ for all $u\in[0,1]^{p+q}$,  hence Corollary \ref{cor: indp theta asymp} holds.

\section{Implementation and Asymptotic Covariance}
\label{sec: computational details}
\noindent \textbf{Estimations of $\hat{\theta}$ and $\hat{m}$}: Since the check loss function in \eqref{objective function} is not differentialable at the origin, the Newton-Raphson algorithm and its variants cannot be directly used to find the solution for linear quantile regression. \cite{koenker2005inequality} proposed several algorithms, such as the interior point algorithm for linear programming, to solve this optimization problem. However, with the layer-by-layer structure of the neural network and the large number of parameters involved, this approach is infeasible for our purpose. 
 We resort to the \textit{Adam} algorithm \citep{kingma2014}, a variant of the \textit{stochastic gradient descent} \citep{robbins1951stochastic}, in the R package \textit{Keras} to solve the optimization problem (\ref{objective function}). This algorithm is widely used in the deep learning field due to its computational and memory efficiency. 
 For our purpose, since we have a parametric and a nonparametric component,  we wrap the linear predictor $\theta^\top X$ and $m(Z)$ together and iteratively estimate  the corresponding parameters simultaneously.  That is, with the neural network $m$ in \eqref{def: DNN2}, we use Adam to update the parameters $\{\theta,W_1,\ldots,W_L\}.$  Here we use the default values in Keras for the initial values $\theta^{(0)}$ and $W_k^{(0)},k=1,\ldots,L$. 
 
The algorithm also requires the specification of  tuning parameters, 
such as the depth $L$, width $\boldsymbol{q}$, step size, minibatch size, the number of iterations and early stopping. Here the minibatch size is defined as the subsample size used to calculate the gradient of the objective function for each iteration, and  early stopping prevents overfitting by specifying the number of iterations to continue
when the model does not improve any more on a hold-out validation dataset. We first hold out 20\% of the training data to select the tuning parameters among a large number of candidates, and then use the selected tuning parameters to redo estimation on the earlier training dataset. Table \ref{tab:hyperparameters} below shows the resulting selected tuning parameters  that are used for the numerical studies in this paper.

 \medskip
\noindent \textbf{Asymptotic Covariance Estimation}: \ To obtain inference for the parameter $\theta_0$, we need to estimate the asymptotic covariance matrix of $\hat{\theta}$ in Theorem \ref{Thm: asymptotic normality} or Corollary \ref{cor: indp theta asymp}. 
For simplicity, we demonstrate how to estimate the asymptotic covariance matrix for the  case  of homoscedastic random errors. 
The first step is to obtain a  density  estimate for $\hat{f}(0)$ from the residuals $\{\hat{\epsilon}_{i}=Y_i-\hat{Y}_i~|~\hat{Y}_i=X_i^\top\hat{\theta}+\hat{m}(Z_i), i=1,\ldots,n\},$ for which we use  \textit{density} in the R package \textit{stats}. 
Then, we employ the deep neural network to estimate the projections ${\varphi}^*_k,k=1,\ldots,p$ empirically, that is,
\begin{equation*}
    \hat{\varphi}^*_k=\mathop{\arg\min}_{\varphi\in\mathcal{M}_1}\frac{1}{n}\sum_{i=1}^{n}\{X_{ik}-\varphi(Z_i)\}^2,
\end{equation*}
where $X_{ik}$ is the $k$-th component of covariates $X_k$ and $\mathcal{M}_1$ is a class of neural networks. Let $\boldsymbol{\hat{\varphi}}^*=(\hat{\varphi}^*_1,\ldots,\hat{\varphi}^*_p)^\top,$ $V_i=X_i-\boldsymbol{\hat{\varphi}}^*(Z_i)$, $\bar{V}=1/n\sum_{i=1}^n V_i$, and
\begin{equation*}
    \hat{\Omega}= \frac{1}{n-1}\sum_{i=1}^{n}(V_i-\bar{V})(V_i-\bar{V})^\top.
\end{equation*}
Finally, we estimate  the asymptotic covariance matrix by 
\begin{equation*}
 \hat{\Sigma} = \frac{\tau(1-\tau)\hat{\Omega}^{-1}}{\hat{f}^2(0)}.
\end{equation*}
For heteroscedastic random errors, we can estimate the corresponding asymptotic covariance matrix by a   bootstrap method, see \cite{feng2011wild} and \cite{wang2018wild} for details.  

\section{Simulations}
\label{sec: simulations}
In this section, we demonstrate the numerical performance of the proposed deep quantile regression method and compare it with  linear quantile regression and partially linear additive quantile regression, abbreviated as LQR and PLAQR, respectively. LQR and PLAQR were implemented with the R packages \textit{quantreg} and \textit{plaqr}, which are publicly available 
at \url{https://cran.r-project.org/package=quantreg} and \url{https://cran.r-project.org/package=plaqr}, respectively. 

\subsection{Simulation I: Homoscedastic Errors}\label{subsec: simulation homo}
We first generated $\tilde{Z}=(\tilde{Z}_1,\ldots, \tilde{Z}_{12})^\top$ from a Gaussian copula on $[0,2]$ with  correlation  parameter $0.5$. Marginally, each coordinate of $\tilde{Z}$ is a uniform distribution on $[0,2]$. We then set $Z=(\tilde{Z}_1,\ldots,\tilde{Z}_{10})^\top$ and $X=(X_1,X_2)^\top$ with $X_1=1(\tilde{Z}_{11}>1)$ and $X_2=\tilde{Z}_{12}$ as covariates. The response $Y$ was generated from 
\begin{equation}
    Y=\theta^\top X + m(Z) +\epsilon,
\end{equation}
where $\theta=(\theta_1,\theta_2)=(1,-1)^\top$, and the error $\epsilon,$ independent of $(X,Z),$ is a Student's t-distribution with zero mean and $3$ degrees of freedom. 
Three choices of  $m$ were implemented:  
\begin{enumerate}[{{\textbf{Case}} 1}]
\item \textbf{(linear)}: $m(z)=0.95\times\sum_{k=1}^{10} z_k$;

\item \textbf{(additive)}: $m(z)=1.1\times\{z_1^3-3z_2^2+2\sin(6\pi z_3) +\log(z_4+0.5) + \sqrt{z_5+2}+e^{z_6/2}+ 0.5(z_7-1+|z_7-1|) +{1}/{(z_8+2)}+2e^{-z_9/2}+\cos(\pi z_{10})\}$;

\item \textbf{(deep)}: $m(z)=0.51\times[z_1z_2+z_2\{1-\cos(\pi z_3z_4) \}+{2\sin(z_5)}/{(|z_5-z_6|+2)}+(z_6+z_7z_8-1)^2+\sqrt{z_9^2+z_{10}^2+2} +\exp\{{\sum_{k=1}^{10}(z_k-1)}/{5}\}].$
\end{enumerate}

The first two cases correspond to, respectively, the LQR and PLAQR model, and the third case is designed for DPLQR.
The factors $0.95, 1.1$ and $0.51$ in each case were scaled to attain a signal-to-noise ratio around $5$. 

For each setting, we generated $Q=160$ 
datasets with respective sample sizes $n=500$ and $2000$ in each dataset. Throughout the simulation, we split the data   into training data and testing data in a 80:20 ratio.  That is, 80\% of the data were used for estimation (including 20\% for tuning) as introduced in Section \ref{sec: computational details},  while 20\%  for evaluating the resulting estimates (test data). 
The performance of $\hat{m}_\tau$ was assessed by the relative mean squared error (RMSE):
\begin{equation}
    RMSE(\hat{m}_\tau)=\frac{\frac{1}{N}\sum_{i=1}^N\{\hat{m}_\tau(Z_i)-m_{\tau}(Z_i)\}^2}{\frac{1}{N}\sum_{i=1}^N\{m_{\tau}(Z_i)\}^2},
\end{equation}
where $\hat{m}_\tau$ and $m_\tau$ are evaluated on the covariates $Z_i,i=1,\ldots,N$ of the test  data. Moreover, with the estimates $\hat{\theta}_\tau$ and $\hat{m}_\tau$, we use $\hat{Y}_i=X_i^\top \hat{\theta}_\tau + \hat{m}_\tau(Z_i)$ to predict  $Y_i$ and  evaluated  its performance through  the mean squared prediction error:
\begin{equation*}
    MSPE(\hat{y})=\frac{1}{N}\sum_{i=1}^N (\hat{Y}_i-Y_i)^2.
\end{equation*}
Here the prediction is also evaluated on the test data.

Table \ref{tab: home theta} presents the biases and standard deviations of the estimates, $\hat{\theta}=(\hat{\theta}_1,\hat{\theta}_2)$, based on 160 simulation runs at three quantile levels $\tau=0.2, 0.5, 0.8$. 
{ In general, both the bias and variance decrease steadily for all three methods as the sample size increases from 500 to 2000.}
As expected, the mean squared error of the resulting estimates are the smallest  at the median ($\tau=0.5$) level.  Under Case 1 (linear) and Case 2 (additive), the proposed DPLOR method performed comparably with the optimal method (LQR and PLAQR respectively) with slightly larger mean squared errors. 
However,  under Case 3 (deep), 
the DPLQR method clearly outperforms LQR and PLAQR. We also construct the $95\%$ confidence intervals for $\theta_1$ and $\theta_2$ based on the estimates of the asymptotic variance in Section \ref{sec: computational details}. 
Table \ref{tab: home cover} reports the empirical coverage probabilities of the 95\% confidence intervals. 
For all three cases, the empirical coverage probabilities of the proposed method generally approach $95\%$ as $n$ increases. Moreover, the proposed method is  comparable to the other two methods under Case 1 (linear) and Case 2 (additive), and has more accurate coverage rates under Case 3 (deep).

The average relative mean squared errors 
of the estimated nonparametric function $\hat{m}$ over $160$ repetitions are given in Table \ref{tab: home m}. They decline with the sample sizes as expected. 
When the true model is Case 3 (deep), the proposed method substantially outperforms LQR and PLAQR, while it performs slightly worse under Case 1 (linear) and Case 2 (additive).

Table \ref{tab: home y} shows the mean of the squared predicted errors 
of the predicted value $\hat{Y}$ based on the median ($\tau=0.5$) regression and  
 reveals that the proposed DPLQR is competitive with the optimal procedure (LQR in  Case 1 and PLAQR in Case 2) and  superior  in Case 3. 

\begin{table}
\caption{\label{tab: home theta}Bias and standard deviation (in parentheses) of $\hat{\theta}$ for the LQR, PLAQR and DPLQR methods under homoscedastic random errors.}\vspace{3mm}
\centering
\fbox{\resizebox{1.0\textwidth}{!}{
\begin{tabular}{ccccccccccc}
&&\multicolumn{3}{c}{$\tau=0.2$}&\multicolumn{3}{c}{$\tau=0.5$}&\multicolumn{3}{c}{$\tau=0.8$} \\
\cmidrule(r){3-5}\cmidrule(r){6-8}\cmidrule(r){9-11}
Case & $n$ & LQR & PLAQR & DPLQR& LQR & PLAQR & DPLQR& LQR & PLAQR & DPLQR \\\hline
\multicolumn{11}{c}{$\theta_1$}\\\hline
{Case 1} &
{500}&0.0068&-0.0053&0.0169&0.0027&-0.0087&0.0109 &-0.0262&-0.0239&-0.0743   \\
(linear)&&(0.1853)&(0.1873)&(0.1848)&(0.1611)&(0.1618)&(0.1602)&(0.2564)&(0.2787)&(0.2806)  \\
 & {2000}&-0.0046&0.0048&0.0145&0.0028&0.0019&0.0160 &0.0208&0.0164&0.0454 \\
 & &(0.0932)&(0.0940)&(0.0927) &(0.0815)&(0.0827)&(0.0810)&(0.1388)&(0.1455)&(0.1442) \\
 {Case 2} 
&{500}&0.0086&0.0079&-0.0182&0.0021&0.0108&-0.0422 &0.0816&0.0088&-0.0831\\
(additive)&&(0.3890)&(0.3309)&(0.3333) &(0.4140)&(0.3686)&(0.3781)&(0.4891)&(0.4744)&(0.4758)  \\
 & {2000}&-0.0272&-0.0088&-0.0012&0.0018&0.0059&0.0017  &0.0523&0.0050&0.0477  \\
 &&(0.2136)&(0.1637)&(0.1680)&(0.2116)&(0.1889)&(0.1925)  &(0.2517)&(0.2343)&(0.2408)  \\

{Case 3} 
&{500} &0.0099&0.0244&0.0180 &-0.0039&-0.0493&0.0066 &-0.4256&-0.0016&-0.0497   \\
(deep)& &(0.2900)&(0.2674)&(0.1919)&(0.2882)&(0.2618)&(0.1839)&(1.2277)&(0.4148)&(0.3176) \\
 & {2000}&-0.0166&-0.0144&0.0057 &0.0106&-0.0160&0.0135 &-0.4222&-0.0434&-0.0056\\
 & &(0.1517)&(0.1414)&(0.0968)&(0.1402)&(0.1382)&(0.0892) &(0.7151)&(0.2260)&(0.1381)
\\\hline
\multicolumn{11}{c}{$\theta_2$}
\\\hline
{Case 1} &
{500}&-0.0135&-0.0236&0.0357&0.0109&0.0191&0.0321&-0.0530&-0.0483&0.2028    \\
(linear)& &(0.1880)&(0.1952)&(0.1977)&(0.1605)&(0.1701)&(0.1634) &(0.3390)&(0.3444)&(0.3976)  \\
 & {2000}&-0.0021&-0.0056&0.0182 &0.0087&0.0093&0.0180&0.0003&-0.0045&0.0952 \\
 & &(0.0882)&(0.0904)&(0.0930)&(0.0749)&(0.0775)&(0.0787)&(0.1641)&(0.1665)&(0.1678)  \\
 {Case 2} 
&{500}&0.0033&0.0021&0.0135&0.0369&-0.0117&0.0258&-0.0140&0.0426&0.1427  \\
(additive)&&(0.3274)&(0.2639)&(0.2813)&(0.4024)&(0.3750)&(0.3785)&(0.6286)&(0.6140)&(0.4951)  \\
 & {2000} &-0.0020&-0.0011&0.0052&-0.0048&-0.036&-0.0042&0.0185&0.0072&0.0121 \\
 &&(0.1800)&(0.1329)&(0.1358)&(0.2096)&(0.2033)&(0.2075)&(0.2782)&(0.2473)&(0.2468)  \\

{Case 3} 
&{500}&-0.0216&-0.0332&0.0387 &-0.0011&-0.0282&0.0802&0.0535&0.0242&0.2552   \\
(deep)&&(0.2886)&(0.2561)&(0.1893)&(0.2955)&(0.2592)&(0.1856)&(1.0797)&(0.3909)&(0.3088) \\
 & {2000}&-0.0035&-0.0058&0.0171  &-0.0064&-0.0127&0.0244 &0.0425&-0.0198&0.0814\\
 & &(0.1523)&(0.1384)&(0.0887)&(0.1695)&(0.1292)&(0.0935)&(0.5007)&(0.2013)&(0.1607)
\\
\end{tabular}
} }
\end{table}

\begin{table}
\caption{\label{tab: home cover} Empirical coverage probability  of the 95\% confidence interval for $\theta=(\theta_1,\theta_2)$ by the LQR, PLAQR and DPLQR    methods under homoscedastic random errors.}\vspace{3mm}
\centering
\fbox{\resizebox{1.0\textwidth}{!}{
\begin{tabular}{ccccccccccc}
&&\multicolumn{3}{c}{$\tau=0.2$}&\multicolumn{3}{c}{$\tau=0.5$}&\multicolumn{3}{c}{$\tau=0.8$} \\
\cmidrule(r){3-5}\cmidrule(r){6-8}\cmidrule(r){9-11}
Case & $n$ & LQR & PLAQR & DPLQR& LQR & PLAQR & DPLQR& LQR & PLAQR & DPLQR \\\hline
\multicolumn{11}{c}{$\theta_1$}\\\hline
{Case 1} &
{500} &0.9500&0.9125&0.9750&0.9188&0.8875&0.9688&0.9188&0.9125&0.9750  \\
 (linear)& {2000}&0.9625&0.9312&0.9688&0.9500&0.9312&0.9625&0.9312&0.9250&0.9625\\
 {Case 2} 
&{500}&0.9000&0.9062&0.9812 &0.9125&0.8750&0.9688 &0.9062&0.8875&0.9688  \\
 (additive)& {2000}&0.8938&0.9250&0.9688&0.9125&0.9375&0.9125 &0.8438&0.9375&0.9625 \\
{Case 3} 
&{500}&0.8750&0.9250&0.9125&0.8875&0.8750&0.9375&0.9688&0.8938&0.8812   \\

 (deep)& {2000}&0.9312&0.8750&0.9375&0.9688&0.8688&0.9562&0.5062&0.8938&0.9438\\
\hline
\multicolumn{11}{c}{$\theta_2$}\\
\hline
{Case 1} &
{500}&0.9312&0.8812&0.9250&0.8938&0.9188&0.9250&0.9062&0.8938&0.8063   \\
(linear) & {2000}&0.9375&0.9312&0.9375&0.9500&0.9188&0.9625&0.9100&0.062&0.9000\\
 {Case 2} 
&{500}&0.8938&0.9375&0.9688&0.8688&0.8875&0.9125&0.8500&0.9312&0.8875   \\
(additive) & {2000} &0.8750&0.9438&0.9562&0.8562&0.9375&0.9125&0.8875&0.9438&0.8938\\
{Case 3} 
&{500} &0.8938&0.9062&0.9125&0.8812&0.8875&0.9125&0.9500&0.8938&0.8250  \\
(deep) & {2000}&0.9062&0.8688&0.9375&0.9688&0.8688&0.9312&0.8000&0.8875&0.9250\\
\end{tabular}
}}
\end{table}

\begin{table}
\caption{\label{tab: home m}Relative mean squared error of $\hat{m}$ for the LQR, PLAQR and DPLQR methods 
under homoscedastic random errors.}\vspace{3mm}
\centering
\fbox{\resizebox{1.0\textwidth}{!}{
\begin{tabular}{ccccccccccc}
&&\multicolumn{3}{c}{$\tau=0.2$}&\multicolumn{3}{c}{$\tau=0.5$}&\multicolumn{3}{c}{$\tau=0.8$} \\
\cmidrule(r){3-5}\cmidrule(r){6-8}\cmidrule(r){9-11}
Case & $n$ & LQR & PLAQR & DPLQR& LQR & PLAQR & DPLQR& LQR & PLAQR & DPLQR \\\hline
{Case 1} &
{500} &0.0019&0.0045&0.0037 &0.0009&0.0020&0.0018&0.0021&0.0037&0.0049   \\
 (linear)& {2000} &0.0004&0.0010&0.0009 &0.0002&0.0004&0.0004&0.0005&0.0008&0.0010\\
 {Case 2} 
&{500}&0.2650&0.1560&0.2361&0.2176&0.1503&0.1962&0.2492&0.2018&0.2199   \\
 (additive)& {2000}&0.2518&0.1362&0.1596&0.1937&0.1175&0.1347&0.2202&0.1647&0.1811 \\
{Case 3} 
&{500}&0.1307&0.1188&0.0796&0.0955&0.0345&0.0183&0.1464&0.0152&0.0132    \\
 (deep)& {2000}&0.1244&0.0996&0.0232&0.0899&0.0258&0.0087 &0.1160&0.0080&0.0053\\
 
\end{tabular}
}}
\end{table}

\begin{table}
\caption{\label{tab: home y} Mean of the squared prediction errors evaluated on the test set for the LQR, PLAQR and DPLQR methods under homoscedastic random errors.}\vspace{3mm}
\centering
\fbox{
\begin{tabular}{ccccc}
Case & $n$ & LQR & PLAQR & DPLQR\\\hline
{Case 1} 
&{500} &3.1028&3.1948&3.1774   \\
 (linear)& {2000} &2.9302&2.9545&2.9484 \\
 {Case 2} 
&{500}&8.1803&6.8469&7.1238  \\

(additive) & {2000} &7.8553&6.1522&6.3471 \\
{Case 3} 
&{500} &6.9982&4.6561&3.8578 \\
(deep) & {2000}&6.3220&3.9613&3.1862\\

\end{tabular}
}
\end{table}

\begin{table}
\caption{\label{tab: hete theta} Bias and standard deviation (in parentheses) of $\hat{\theta}$ for the LQR, PLAQR and DPLQR methods under heteroscedastic random errors.}\vspace{3mm}
\centering
\fbox{\resizebox{\textwidth}{!}{
\begin{tabular}{ccccccccccc}
&&\multicolumn{3}{c}{$\tau=0.2$}&\multicolumn{3}{c}{$\tau=0.5$}&\multicolumn{3}{c}{$\tau=0.8$} \\
\cmidrule(r){3-5}\cmidrule(r){6-8}\cmidrule(r){9-11}
Case & $n$ & LQR & PLAQR & DPLQR& LQR & PLAQR & DPLQR& LQR & PLAQR & DPLQR \\\hline
\multicolumn{11}{c}{$\theta_1$}\\\hline
{Case 4} 
&{500}&0.0351&0.0459&0.0525 &-0.0127&0.0037&0.0051&0.0693&0.0960&-0.1815   \\
(linear)&&(0.4054)&(0.4443)&(0.4042)&(0.3070)&(0.3240)&(0.3102)&(0.5661)&(0.6115)&(0.5699)    \\
 &{2000}&-0.0332&-0.0436&0.0085 &-0.0166&-0.0089&0.0130&0.0005&0.0040&0.0294 \\
 &&(0.2254)&(0.2414)&(0.2249)&(0.1736)&(0.1837)&(0.1779)&(0.2866)&(0.3095)&(0.2908)  \\
{Case 5} 
&{500} &0.0437&0.0785&-0.0302 &0.0058&0.0203&-0.0671 &0.0855&-0.1677&-0.2673   \\
(additive)& &(0.5168)&(0.4704)&(0.4799) &(0.5006)&(0.4305)&(0.4435)&(0.7398)&(0.7104)&(0.7271) \\
 & {2000}&0.1465&0.1086&0.1424 &-0.0041&-0.0180&0.0027&0.1163&-0.1208&-0.0729 \\
 &&(0.2750)&(0.2437)&(0.2452)&(0.2547)&(0.2253)&(0.2292) &(0.3699)&(0.3562)&(0.3625)   \\

{Case 6} 
&{500}&0.0511&0.0346&0.0494&0.0784&-0.0472&-0.0160&1.4354&-0.0252&-0.0695\\
(deep)& &(0.5284)&(0.4786)&(0.3563)&(0.4040)&(0.3829)&(0.2911)&(1.4179)&(0.6577)&(0.4847)  \\
 & {2000}&0.0629&0.0484&0.0230&0.0910&-0.0203&0.0054 &1.6135&-0.0159&-0.0022\\
 &&(0.2613)&(0.2353)&(0.1842)&(0.2188)&(0.1872)&(0.1354)&(0.6215)&(0.3223)&(0.2313) \\\hline
\multicolumn{11}{c}{$\theta_2$}
\\\hline
{Case 4} 
 & {500}&-0.0507&-0.0483&0.1967&-0.0111&-0.0066&0.2677&-0.0094&0.0116&0.3887   \\
 (linear)&&(0.3678)&(0.3695)&(0.3688)&(0.3234)&(0.3300)&(0.3346)&(0.6264)&(0.6352)&(0.6381)   \\
 & {2000} &0.0101&0.0097&0.0663 &0.0071&0.0048&0.0607 &-0.0538&-0.0400&0.1558\\
 &&(0.1846)&(0.1875)&(0.1861)&(0.1541)&(0.1617)&(0.1592)&(0.2958)&(0.2919)&(0.3009)\\
{Case 5} 
&{500} &0.0660&0.0851&0.1920 &-0.0376&-0.0237&0.0064 &-0.1001&-0.1008&0.2185   \\
(additive)&&(0.4954)&(0.4484)&(0.4570)&(0.5132)&(0.4182)&(0.4200)&(0.8267)&(0.7533)&(0.8010)   \\
 & {2000}&0.0680&0.0712&0.0854 &-0.0034&-0.0060&0.0087 &-0.1277&-0.1116&-0.0383\\
 &&(0.2565)&(0.2150)&(0.2186)&(0.2755)&(0.2083)&(0.2125)&(0.3939)&(0.3686)&(0.3767)  \\

{Case 6} 
&{500}&0.0645&0.0987&0.0772&0.0697&0.0156&0.0865 &0.0120&-0.0711&0.2220   \\
(deep)& &(0.4301)&(0.3816)&(0.2843)&(0.3965)&(0.3673)&(0.2504)&(0.9971)&(0.5929)&(0.3597) \\
 & {2000}&0.0414&0.0302&0.0285&0.0489&-0.0179&0.0207 &-0.0721&-0.0335&0.0931\\
 &&(0.2253)&(0.1993)&(0.1391)&(0.1917)&(0.1838)&(0.1187)&(0.4851)&(0.3062)&(0.1686)\\
\end{tabular}
} }
\end{table}

\begin{table}
\caption{\label{tab: hete cover} Empirical coverage probability of the 95\% confidence interval for $\theta=(\theta_1,\theta_2)$ by the LQR, PLAQR and DPLQR methods under heteroscedastic random errors.}\vspace{3mm}
\centering
\fbox{\resizebox{\textwidth}{!}{
\begin{tabular}{ccccccccccc}
&&\multicolumn{3}{c}{$\tau=0.2$}&\multicolumn{3}{c}{$\tau=0.5$}&\multicolumn{3}{c}{$\tau=0.8$} \\
\cmidrule(r){3-5}\cmidrule(r){6-8}\cmidrule(r){9-11}
Case & $n$ & LQR & PLAQR & DPLQR& LQR & PLAQR & DPLQR& LQR & PLAQR & DPLQR \\\hline
\multicolumn{11}{c}{$\theta_1$}\\
\hline
{Case 4} 
&{500}&0.9188&0.9188&0.9125 &0.9125&0.9562&0.9750&0.9125&0.8875&0.8688  \\

(linear) & {2000}&0.9375&0.9188&0.9250&0.9562&0.9250&0.9188&0.9375&0.8938&0.9125\\

{Case 5} 
&{500}&0.8878&0.9062&0.9625&0.8812&0.8878&0.8812&0.8625&0.9250&0.9750   \\

(additive) & {2000}&0.9125&0.9375&0.9250&0.9062&0.9250&0.9125&0.8250&0.9562&0.9375\\

{Case 6} 
&{500} &0.8875&0.8750&0.8938 &0.9000&0.8625&0.9188&0.9062&0.8625&0.9062    \\

(deep) & {2000} &0.9063&0.8938&0.9250&0.8875&0.8875&0.9188&0.8750&0.9188&0.9250 \\
 \hline
\multicolumn{11}{c}{$\theta_2$}\\
\hline
{Case 4} 
&{500} &0.8812&0.8750&0.8938  &0.8688&0.8812&0.8625&0.8812&0.9125&0.8750 \\

(linear) & {2000} &0.9188&0.9062&0.9062  &0.9375&0.9125&0.9062&0.9312&0.9250&0.9125 \\

{Case 5} 
&{500}&0.8750&0.9125&0.9062&0.9000&0.9125&0.8750&0.8875&0.8875&0.8750  \\

(additive) & {2000}&0.8312&0.9125&0.9125&0.8620&0.9312&0.9250&0.9000&0.9125&0.9062\\

{Case 6} 
&{500} &0.9062&0.8875&0.9125&0.8750&0.9000&0.9000&0.9250&0.9250&0.8812   \\

(deep) & {2000}&0.8875&0.8688&0.9188&0.8625&0.8688&0.9250&0.8500&0.8562&0.9125 \\
\end{tabular}
}}
\end{table}

\begin{table}
\caption{\label{tab: hete m} Relative mean squared error of $\hat{m}$ for the LQR, PLAQR and DPLQR methods under heteroscedastic random errors.}\vspace{3mm}
\centering
\fbox{\resizebox{\textwidth}{!}{
\begin{tabular}{ccccccccccc}
&&\multicolumn{3}{c}{$\tau=0.2$}&\multicolumn{3}{c}{$\tau=0.5$}&\multicolumn{3}{c}{$\tau=0.8$} \\
\cmidrule(r){3-5}\cmidrule(r){6-8}\cmidrule(r){9-11}
Case & $n$ & LQR & PLAQR & DPLQR& LQR & PLAQR & DPLQR& LQR & PLAQR & DPLQR \\\hline
{Case 4} 
&{500}&0.0108&0.0276&0.0199&0.0032&0.0089&0.0062&0.0058&0.0145&0.0103\\

(linear) & {2000}&0.0026&0.0063&0.0050&0.0009&0.0021&0.0016&0.0012&0.0029&0.0026\\

{Case 5} 
&{500} &0.2799&0.2075&0.2241&0.2306&0.1796&0.2034&0.2679&0.2206&0.2296\\

 (additive)& {2000}&0.2514&0.1432&0.1716 &0.1948&0.1232&0.1440&0.1791&0.1389&0.1492\\

{Case 6} 
 &{500}&0.2830&0.2362&0.1895&0.1086&0.0685&0.0329&0.1098&0.0308&0.0187\\

(deep)&{2000}&0.1876&0.1450&0.0697&0.0908&0.0388&0.0186&0.0865&0.0157&0.0072\\

\end{tabular}
}
}
\end{table}

\begin{table}
\caption{\label{tab: hete y}Mean of the squared prediction errors evaluated on the test set for the LQR, PLAQR and DPLQR methods under heteroscedastic random errors.}\vspace{3mm}
\centering
\fbox{\resizebox{0.5\textwidth}{!}{
\begin{tabular}{ccccc}
Case & $n$ & LQR & PLAQR & DPLQR\\\hline
{Case 4} 
&{500}&18.0770&18.7186&18.3779  \\

(linear) & {2000} &17.8294&17.9436&17.9470 \\
{Case 5} 
&{500} &21.8791&19.9281&20.5156\\
(additive) & {2000}&21.6721&19.8431&20.0425\\

{Case 6} 
&{500} &18.5553&17.0905&16.0694\\
(deep) & {2000} &18.4154&16.8250&15.7411 \\

\end{tabular}
}}
\end{table}

\subsection{Simulation II: Heteroscedastic Errors}
\label{subsec: simulation hete}
We also studied the performance of the proposed method for heteroscedastic errors. 
The covariates $U=(X,Z)$, coefficient $\theta$ and nonparametric function $m$  are similar to the settings in Section 5.1 but the response $Y$ now comes from the regression model:
\begin{equation*}
    Y = X^\top\theta + m(Z) + \sigma_1(X,Z)\epsilon.
\end{equation*}
Here $\epsilon$ follows the Student’s t-distribution with zero mean and 3 degrees of freedom. The function $\sigma_1(X,Z)$ has 
the following three settings: 
\begin{enumerate}[{{\textbf{Case}}}]
\item \textbf{4 (linear)}: $\sigma_1(x,z)=(x_1+x_1+\sum_{k=1}^{10}z_k)/5;$
\item \textbf{5 (additive)}: $\sigma_1(x,z)=(x_1+x_1+\sum_{k=1}^{10}|z_k-0.2|)/3.6;$
\end{enumerate}
\begin{enumerate}[{{\textbf{Case 6 (deep)}:}}]
\item $\sigma_1(x,z)= (x_1+x_1)/3+3\Phi\big({\sum_{k=1}^{10}(z_k-1) }/{5}\big)$ with  the cumulative distribution function $\Phi(\cdot)$ of the standard normal distribution.
\end{enumerate}
These lead to 
$\theta_\tau=\theta+t_\tau\theta^*$ and 
$m_\tau(z)=m(z)+t_\tau m^*(z)$ with $t_\tau$ being the $\tau$ quantile of  Student’s t-distribution with zero mean and degree of freedom 3. 
The simulation results, which are summarized in Table \ref{tab: hete theta} - Table \ref{tab: hete y}, are  comparable  to those in  Simulation I in Section \ref{subsec: simulation homo}. 

In summary, when the true model is  linear or partially linear additive quantile regression, our method is  competitive 
for  both the parametric coefficients and nonparametric function estimates, and the coverage probabilities for the parametric coefficients are close to the 95\% nominal level as sample sizes increase. Furthermore, the proposed method is superior to the LQR and PLAQR methods when the true model comes from the deep partially linear quantile regression. 

\section{Applications to Real Data}
\label{sec: application}

\subsection{Concrete Compressive Strength Data}
We apply the proposed methodology, along with the competing methods, to the Concrete Compressive Strength Data Set \citep{yeh1998modeling} available on the UCI machine learning repository. The data consist of $n=1030$ observations, with the response being a continuous variable of concrete compressive strength (CCS), and eight covariates: $Z_1$(cement), $Z_2$(water), $Z_3$(fly ash), $Z_4$(blast furnace slag), $Z_5$(superplasticizer), $Z_6$(coarse aggregate), $Z_7$(fine aggregate) and $Z_8$(age of the mixture in days), of which the first seventh covariates are the ingredients in high-performance concrete (HPC). For concrete technology, the water-cement ratios (WCR) has been recognized as the most useful and significant advancement for CCS \citep{yeh1998modeling}. Here we not only explore the association between WCR and CCS, but also predict the CCS of HPC from the covariates.

As in the simulations, we  model the data with four approaches: \textbf{(a)} the proposed deep partially linear quantile regression (DPLQR), \textbf{(b)} linear quantile regression (LQR), \textbf{(c)} partially 
linear additive quantile  regression (PLAQR), and \textbf{(d)} deep nonparametric quantile regression ($\text{DNQR}$, see \cite{jantre2020quantile,padilla2020quantile}). 
Note that model (d) does not offer direct treatment effects.

For the data, we treat log CCS as the response,  WCR, i.e, $Z_2/Z_1$, as the linear predictors $X$, and $(Z_3,Z_4,Z_5,Z_6,Z_7,Z_8)^\top$ as the predictors $Z$ in models (a), (b) and (c). In model (d), we nonparametriclly regress response log CCS on all covariates $(Z_1,\ldots,Z_8)$  and implement it via deep learning. We use 80\% of the data to train and tune the model and hold out the rest 20\% of data to assess the prediction performance of the four methods. 

Table \ref{tab: est&pred} shows the numerical results at the median level ($\tau=0.5$). The 95\% confidence interval for $\theta$ of each method suggests that there is a strong association between WCR and CCS. The negative estimates  further support the  Abrams rule in civil engineering that  increase in the WCR tends to decrease the strength of concrete \citep{gorse2012dictionary}. 
Among the first three approaches, DPLQR produced a shorter 95\% confidence interval for $\theta$ than the LQR and PLAQR methods.  The prediction results in Table \ref{tab: est&pred} further reveal that our method not only improves the prediction accuracy  substantially but also has the smallest standard deviation.
Although the proposed model is a submodel of the DNQR in (d), its performance in prediction is comparable. Thus,  compared to the fully nonparametric approach (d), the partially linear approach (c) trades a small amount of prediction accuracy for interpretibility and has the best performance among the three interpretable approaches (a)-(c). 



\begin{table}
\caption{\label{tab: est&pred} Estimation  and prediction results for the Concrete Compressive Strength Data. CI: confidence interval; MEAN: mean of the squared prediction errors; SD: standard deviation of the squared prediction errors.}\vspace{3mm}

\centering
\fbox{
\begin{tabular}{ccccc}
&\multicolumn{2}{c}{Estimation}&\multicolumn{2}{c}{Prediction error} \\
\cmidrule(r){2-3}\cmidrule(r){4-5}
&$\hat{\theta}$&95\% CI&MEAN&SD\\
\hline
LQR&-1.3043&[-1.4803, -1.2147]& 0.1718&0.2557\\
PLAQR&-1.0959&[-1.2718, -1.0241]&0.0708&0.1080 \\
DPLQR&-1.2627&[-1.3248, -1.2006]&0.0319&0.0740\\
{DNQR}&-&-&0.0308&0.0615\\
 
\end{tabular}
}%
\end{table}

\subsection{Boston Housing Data}
The Boston Housing Data, available from the R package \textit{mlbench}, is a  benchmark dataset for quantile regression analysis. In \cite{harrison1978hedonic}, 506 observations were examined to study 
the housing prices based on various demographic and socioeconomic predictors. The  variables are: $Y$ (median hoouse price), $X_1$ (per capita crime rate by town), $X_2$ (a river boundary indicator), $X_3$ (proportion of non-retail business acres per town), $X_4$ (proportion of residential land zoned for lots), $X_5$ (nitrogen oxides concentration), $X_6$ (average number of rooms per dwelling), $X_7$ (proportion of owner-occupied units built prior to 1940), $X_8$ (weighted mean of distances to five Boston employment centers), $X_9$ (index of accessibility to radial highways), $X_{10}$ (full-value property-tax rate), $X_{11}$ (pupil-teacher ratio by town), $X_{12}$ (the proportion of black individuals by town), $X_{13}$ (the percentage of the population classified as lower status). 

To study the effect of the crime rate on house price, we choose $X_1$ (per capita crime rate by town) and the binary coavriate $X_2$ (a river boundary indicator) as the vector predictors and all other {continuous} covariates as the nonparametric predictors. 

\begin{equation*}
    \log Y = \theta_1 X_1 +\theta_2 X_2 + m(X_3,\ldots,X_{13}) +\epsilon.
\end{equation*}
Here the function $m$ is modelled as a linear, additive and nonparametric function, which corresponds to the LQR, PLAQR and DPLQR models, respectively. We also include  the DNQR method, which treats all thirteen covariates as components of a nonparametric regression model, i.e. $\log Y= m(X_1, X_2, X_3, \ldots,X_{13}) +\epsilon.$ 

The estimates $\hat{\theta}_1$ from the median regression are summarized in Table \ref{tab: est&pred boston},  revealing that the crime rate has a significant  effect on the price of a house and 
house prices are higher in areas with lower  crime rates. Table \ref{tab: est&pred boston} also displays the mean and standard deviation of squared prediction errors ( hold-out 20\% as test set). 
The proposed method is considerably better than the LQR and PLADR methods. Furthermore, we note that the DNQR method leads to a larger  squared prediction error than  our method.

\begin{table}
\caption{\label{tab: est&pred boston} Estimation and prediction results for the Boston Housing Data. CI: confidence interval; MSPE: mean of the squared prediction errors; SD: standard deviation of the squared prediction errors.}\vspace{3mm}
\centering
\fbox{
\begin{tabular}{lcccc}
&\multicolumn{2}{c}{Estimation}&\multicolumn{2}{c}{Prediction} \\
\cmidrule(r){2-3}\cmidrule(r){4-5}
&$\hat{\theta}_1$&95\% CI&MEAN&SD\\
\hline
LQ&-0.0093&[-0.0274, -0.0081]&0.0799&0.2048\\
PLAQR&-0.0112&[-0.0309,  -0.0082]&0.0529&0.1314 \\
DPLQR&-0.0117&[-0.0137, -0.0096]&0.0272&0.0559\\
$\text{DNQR}$&-&-&0.0283&0.0546\\
\end{tabular}
}
\end{table}

\begin{table}
\caption{\label{tab:hyperparameters} Tuning parameters in the simulations and data applications}\vspace{3mm}
\centering
\fbox{\resizebox{\textwidth}{!}{
\begin{tabular}{lcccccccc}
\hline
&\multicolumn{2}{c}{Case 1\&4 (linear)}&\multicolumn{2}{c}{Case 2\&5 (additive)}&\multicolumn{2}{c}{Case 3\&6 (deep)}&\multirow{2}{*}{Concrete Data}&\multirow{2}{*}{Housing Data} \\
\cmidrule(r){2-3}\cmidrule(r){4-5}\cmidrule(r){6-7}
&500&2000&500&2000&500&2000&\\\hline
Depth&2&3&3&3&2&3&3&3\\
Width&16&32&10&20&20&32&32&32\\
Epoch&500&500&500&500&600&600&1000&500\\
Minibatch&64&64&64&64&128&128&64&64\\
Early stop&50&50&50&50&100&100&100&50\\
Learning rate&0.01/0.02&0.01/0.02&0.009/0.01&0.009/0.02&0.01/0.02&0.01/0.02&0.009&0.02\\

\end{tabular}
}}
\end{table}

\section{Conclusion}
\label{sec: funture works} 
We provide an interpretable-yet-flexible deep learning model with partially linear quantile regression, where we leverage the neural networks to represent the nonparametric function and the linear predictor to obtain inference. The proposed method is able to detect the parsimonious structure of the data automatically, thereby producing a better convergence rate for the nonparametric estimator $\hat{m}$ than conventional nonparametric smoothing methods. Furthermore, the estimator of the parameter $\theta_0$ {attains} $\sqrt{n}$-consistency and asymptotic normality. These substantially distinguish our method from  neural networks for  nonparametric regression \citep{padilla2020quantile,schmidt2020nonparametric}, and also open up a myriad of research opportunities for semiparametric regression models. 

A possible extension  is to investigate the quantile regression process instead of fitting  a  quantile level $\tau$. \cite{chao2017quantile} and \cite{belloni2019conditional} studied convergence results uniformly on $\tau$ for quantile functions approximated by linear combinations of basis functions obtained, e.g. from  polynomial, Fourier, spline and wavelet bases.  However, their approaches cannot easily be extended to the deep learning setting because of the  layer structure in a neural network. To further investigate this therefore  will be an interesting future project. 

As we focus in this paper on a fixed but moderate size of the linear covariates $X$,  future work of interest is to study DPLQR with high-dimensional covariates, where the number of linear covariates may grow at a certain rate with  sample size. A special case for PLAQR was studied in \cite{sherwood2016partially}, which may shed some light on extending  the DPLQR approach.  

\section{Proofs of Theorems}
\label{sec: proofofthm}
\noindent \textbf{Proof of Theorem \ref{Thm: convergency rate}.}
Let  $\hat{\beta} =(\hat{\theta},\hat{m})$, $\beta_0=(\theta_0,m_0)$ and $d(\beta_1,\beta_2)=[\mathbb{E}\{x^\top\theta_1+m_1(Z)-x^\top\theta_2-m_2(Z)\}^2]^{1/2}$, for any $\beta_1=(\theta_1,m_1)$ and $\beta_2=(\theta_2,m_2)$. We first show that
\begin{equation*}
    d(\hat{\beta},\beta_0)\xrightarrow{p} 0, ~\text{as}~ n\rightarrow\infty.
\end{equation*}

Choose some large $C>0,$ such that $\|m_0\|_{L^2([0,1]^{q})}<C$ and $\|\theta_0\|<C$ with $\|\cdot\| $ being the Euclidean norm of a vector. Let $\mathbb{R}_{C}^{p}=\{\theta\in\mathbb{R}^{p}~|~\|\theta\|<C \}$ and $\mathcal{M}_{C}=\mathcal{M}(s,L,\boldsymbol{q},C)$ in (\ref{def: sparse NN}). Define
\begin{equation}\label{minimizerC}
    \hat{\beta}_{C}=\mathop{\arg\min}_{\beta\in \mathbb{R}_{C}^{p}\times\mathcal{M}_C} L_{ n}(\beta),
\end{equation}
where $L_{ n}(\beta)=1/n\sum_{i=1}^{n}\rho_{\tau}(Y_i-X_i^\top\theta-m(Z_i))$ for $\beta=(\theta,m)$ and $\rho_\tau$ defined in (\ref{objective function}).
It is easy to show, by contradiction, that $\mathbb{P}(\|\hat{\theta}\|<C,\|\hat{m}\|_{L^2([0,1]^q)}<C)\rightarrow 1,$ as $C\rightarrow\infty$. 
Thus it suffices to verify that $\hat{\beta}_{C}$ is consistent for large enough $C>0$, i.e., $d(\hat{\beta}_{C},\beta_0)\xrightarrow{p} 0,$ as $n\rightarrow\infty.$ 

By Lemma 5 in \cite{schmidt2020nonparametric} and the fact $|\rho_\tau(u)-\rho_\tau(v)|\le 2|u-v|$ for all $u,v\in\mathbb{R}$, we know that $\{\rho_{\tau}(Y-\theta^\top X - m(Z))~|~\beta=(\theta,m)\in\mathbb{R}^{p}_C\times\mathcal{M}_C\}$ is  \textit{P-Glivenko-Cantelli}. 
Hence
\begin{equation}\label{pfthm1: uniform consistemcy}
    \mathop{\sup}_{\beta \in \mathbb{R}^{p}_{C}\times\mathcal{M}_C}|L_{n}(\beta)-L_{0}(\beta)|\xrightarrow{p}0,~\text{as}~n\rightarrow \infty,
\end{equation}
where $L_{0}(\beta)=\mathbb{E}\rho_{\tau}(Y-X^\top\theta-m(Z)).$

By Assumption (A\ref{assump: covariates}) and (A\ref{assump: error pdf}), a similar proof for equation (C.44) in \cite{belloni2019conditional} implies that, for any $\epsilon>0,$ 
\begin{equation}\label{thm15}
    \mathop{\inf}_{ \substack{d(\beta,\beta_0)>\epsilon,\\\beta\in\mathbb{R}^{p}_{C}\times\mathcal{M}_C } } L_{0}(\beta)>L_{0}(\beta_0).
\end{equation}
For the true function $m_0,$ let
\begin{equation*}
    m^* = \mathop{\arg\min}_{m\in\mathcal{M}_C}\|m-m_0\|_{L^2([0,1]^q)} ~\text{and}~\beta^*=(\theta_0,m^*).
\end{equation*}
Then, we have $L_{n}(\hat{\beta}_{C})\le L_n(\beta^*)$ by the definition of $\hat{\beta}_{C}$ in \eqref{minimizerC}. 
This and \eqref{pfthm1: uniform consistemcy}, \eqref{thm15} imply that
\begin{equation*}
    d(\hat{\beta}_{C},\beta^*)\rightarrow 0,~ \text{as} ~n\rightarrow\infty.
\end{equation*}
On the other hand, by equation (26) in \cite{schmidt2020nonparametric}, we have
\begin{equation}\label{pf: apperror}
    \|m^*-m_0\|_{L^2([0,1]^q)}=O(r_n).
\end{equation}
It follows that
\begin{equation*}
    d(\hat{\beta}_{C},\beta_0)\le d(\hat{\beta}_{C},\beta^*) + \|m^*-m_0\|_{L^2([0,1]^q)}\rightarrow 0,~ \text{as} ~n\rightarrow\infty.
\end{equation*}
This completes the proof of the consistency of $\hat{\beta}$.

Next we prove $d(\hat{\beta},\beta_0)=O_{p}(r_n\log^2 n). $ 
Write $R= L\prod_{k=0}^{L}(q_{k}+1)\sum_{k=1}^{L}q_{k-1}q_{k}$ and
\begin{equation}\label{A_delta}
\mathcal{A}_{\delta} = \{\beta \in \mathbb{R}^{p}_{C}\times \mathcal{M}(  s,L,\textbf{q}, C)~|~ \delta/2\le d(\beta, \beta^*)\le \delta \}.
\end{equation}
 We verify that, for any $\delta>0$,
\begin{equation}\label{proof: donsker}
   \mathbb{E}^*[ \mathop{\sup}_{\beta\in\mathcal{A}_{\delta}}\sqrt{n}\{(L_n-L_0)(\beta^*)-(L_n-L_0)(\beta)\}]\lesssim\phi_n(\delta),
\end{equation}
where $\mathbb{E}^*$ is an outer measure, $\phi_n(\delta)=\delta\sqrt{s\log\frac{R}{\delta}} + \frac{s}{\sqrt{n}}\log \frac{R}{\delta}$, and $a_n\lesssim b_n$ means $a_n\le c b_n$ for some constant $c>0$.

 Denote $\rho_{\tau}(\beta)=\rho_{\tau}(Y-X^\top\theta-m(Z))$ and $\mathcal{B}_\delta=\{\rho_{\tau}(\beta^*)-\rho_{\tau}(\beta)~|~\beta\in\mathcal{A}_{\delta}\}.$ 
 For any $\beta,\beta_1\in\mathcal{A}_{\delta}$, we have 
 $\mathbb{E}|\rho_\tau(\beta)-\rho_\tau(\beta_1)|^2\le 4d^2(\beta,\beta_1)$.  
 Lemma 5 in \cite{schmidt2020nonparametric} then implies that
\begin{equation*}
  \log(1+ \mathcal{N}_{[~]}(\epsilon,\mathcal{B}_{\delta},L^{2}(P)))\lesssim s\log\frac{R}{\epsilon},
\end{equation*}
where $\mathcal{N}_{[~]}(\epsilon,\mathcal{B}_{\delta},L^{2}(P))$ is the bracket number of $\mathcal{B}_{\delta}$ with $L^{2}(P)$ norm. It follows that
\begin{equation*}
    J_{[~]}(\delta,\mathcal{B}_{\delta})=\int_{0}^{\delta}\sqrt{1+ \mathcal{N}_{[~]}(\epsilon,\mathcal{B}_{\delta},L^{2}(P))}d\epsilon\lesssim \delta\sqrt{s\log\frac{R}{\delta}}.
\end{equation*}
By Lemma 3.4.2 of \cite{van1996weak}, we conclude that 
\begin{equation*}
\begin{aligned}
   \mathbb{E}^* [\mathop{\sup}_{\beta\in\mathcal{A}_{\delta}}\sqrt{n}\{(L_n-L_0)(\beta^*)-(L_n-L_0)(\beta)\}]&= \mathbb{E}^*[ \mathop{\sup}_{\beta\in\mathcal{A}_{\delta}}\sqrt{n}(\mathbb{P}_n-\mathbb{P})\{\rho_{\tau}(\beta^*)-\rho_{\tau}(\beta)\}]\\
   &\lesssim J_{[~]}(\delta,\mathcal{B}_{\delta})\Big\{ \frac{J_{[~]}(\delta,\mathcal{B}_{\delta})}{\delta^2\sqrt{n}}+1\Big\}\\
   &=\phi_n(\delta).
   \end{aligned}
\end{equation*}

Let $\eta_n=r_n\log^2n.$ It is clear that 
\begin{equation}\label{convex_sqrt}
    \frac{1}{\eta_n^2}\phi_n(\eta_n)\lesssim\sqrt{n}~\text{and}~L_n(\hat{\beta}_{C})\le L_{n}(\beta^*).
\end{equation}
Then {with \eqref{proof: donsker}, \eqref{convex_sqrt} and} Theorem 3.4.1 of \cite{van1996weak}, we have $d(\hat{\beta}_C,\beta^*)=O_p(\eta_n).$ Hence, It follows from \eqref{pf: apperror} that  $d(\hat{\beta},\beta_0)=O_p(r_n\log^2 n).$ 

Moreover, by Assumption (A\ref{assump: error pdf}) and the definition of $\boldsymbol{\varphi}^*$, we have
\begin{equation*}
\begin{aligned}
    d^2(\hat{\beta},\beta_0)&=\mathbb{E}\{X^\top(\hat{\theta}-\theta_0) + \hat{m}(Z)-m_0(Z)\}^2\\
    &\ge \frac{1}{c_0}\mathbb{E}[f(0|U)\{X^\top (\hat{\theta}-\theta_0) + \hat{m}(Z)-m_0(Z)\}^2]\\
    &=\frac{1}{c_0}\mathbb{E}[f(0|U) \{(X-\boldsymbol{\varphi}^*(Z))^\top (\hat{\theta}-\theta_0) + (\hat{\theta}-\theta_0)^\top \boldsymbol{\varphi}^*(Z)+ \hat{m}(Z)-m_0(Z)\}^2]\\
    &=\frac{1}{c_0}\mathbb{E}[f(0|U)\{ (X-\boldsymbol{\varphi}^*(Z))^\top(\hat{\theta}-\theta_0)\}^2]\\
    &~~~+\frac{1}{c_0}\mathbb{E}[f(0|U)\{(\hat{\theta}-\theta_0)^\top \boldsymbol{\varphi}^*(Z) + \hat{m}(Z)-m_0(Z)\}^2].
\end{aligned}
\end{equation*}
Since the matrix $\mathbb{E}[f(0|U)\{X-\boldsymbol{\varphi}^*(Z)\}\{X-\boldsymbol{\varphi}^*(Z)\}^\top]$ is positive definite, it follows that $\|\hat{\theta}-\theta_0\|=O_p(r_n\log^2 n)$ and thus $\|\hat{m}-m_0\|_{L^2([0,1]^q)}=O_p(r_n\log^2 n)$. This completes the proof.\\ 



\noindent \textbf{Proof of Theorem \ref{Thm: minimax lower bound}.}

For simplicity, we only consider the proof for the median quantile regression case, when $\tau=0.5$. To derive the minimax lower bound, it suffices to show that, when the error $\epsilon$ is  the standard normal distribution and the parameter $\theta_0$ is known and fixed, {there exists a subset} $\mathcal{H}^{*}$ of $\mathcal{H}(J,\boldsymbol{\gamma},{\boldsymbol{d}},\boldsymbol{\bar{d}},B)$ in Assumption (A\ref{assump: param space}), {such that}
\begin{equation}\label{pf: subset minimax}
    \mathop{\lim}_{C\rightarrow\infty}\mathop{\lim}_{n\rightarrow\infty}\mathop{\inf}_{\hat{m}}\mathop{\sup}_{m_0\in\mathcal{H}^*}\mathbb{P}_{(\theta_0,m_0,f)}\big(\|\hat{m}-m_0\|_{L^2([0,1]^q)}\ge C r_n\big)=1,
\end{equation}
where $f$ is the probability density function of standard normal distribution.

Let $KL(\cdot,\cdot)$ be the Kullback-Leibler distance. Suppose that there exists $m^{(0)},\ldots,m^{(N)}\in \mathcal{H}^*$ with $N$ increasing with $n,$ such that for some constants $c_1,c_2>0,$
\begin{equation}\label{pf: m_jk seperate}
    \|m^{(j)}-m^{(k)}\|_{L^2([0,1]^q)}\ge 2c_1 r_n ~\text{for any}~ 0\le j<k\le N, 
\end{equation}
and
\begin{equation*}
    \frac{1}{N}\sum_{i=1}^{N}KL(P_j,P_0)\le c_2\log N,
\end{equation*}
where $P_j$ is the laws corresponding to $(\theta_0,m^{(j)},f)$ for $j=0,\ldots,N,$ respectively. 
Then,  Theorem 2.5 of \cite{tsybakov2008introduction} implies that
\begin{equation*}
    \mathop{\inf}_{\hat{m}}\mathop{\sup}_{m_0\in\mathcal{H}^*}\mathbb{P}_{(\theta_0,m_0,f)}\big(\|\hat{m}-m_0\|_{L^2([0,1]^q)}\ge c_1 r_n\big)\ge \frac{\sqrt{N}}{1+\sqrt{N}}\Big(1-2c_2-\sqrt{\frac{2c_2}{\log N}}\Big).
\end{equation*}
The result \eqref{pf: subset minimax} thus follows.

Note that the likelihood function of $P_j$ with the data  $\{(Y_i,X_i,Z_i)~|~i=1,\ldots,n\}$ and $m^{(j)}$ satisfies
\begin{equation*}
    P_j = \prod_{i=1}^{n} \{f(Y_i-\theta_0^\top X_i-m^{(j)}(Z_i))g(X_i,Z_i)\},
\end{equation*}
where $g$ is the joint probability density of $(X,Z).$ It follows that, if the density of $Z$ is uniformly bounded by a constant $c_3>0,$ then
\begin{equation*}
    \begin{aligned}
       \frac{1}{N}\sum_{i=1}^{N}KL(P_j,P_0)&=\frac{n}{2N}\sum_{j=1}^{N}\mathbb{E}\{m^{(j)}(Z)-m^{(0)}(Z) \}^2\\&\le \frac{c_3n}{2N}\sum_{j=1}^{N}\|m^{(j)}-m^{(0)}\|_{L^2([0,1]^q)}^2.
    \end{aligned}
\end{equation*}
Then by a similar construction as in the proof of Theorem 3 of \cite{schmidt2020nonparametric}, 
there exist $m^{(0)},\ldots,m^{(N)}\in \mathcal{H}(J,\boldsymbol{\gamma},{\boldsymbol{d}},\boldsymbol{\bar{d}},B)$ and constants $c_1,c_2>0$ satisfying both \eqref{pf: m_jk seperate}
and
\begin{equation*}
    \frac{n}{N}\sum_{j=1}^{N}\|m^{(j)}-m^{(0)}\|_{L^2([0,1]^q)}^2\le c_2\log N.
\end{equation*}
The proof is thus complete.\\

\noindent \textbf{Proof of Theorem \ref{Thm: asymptotic normality}.}
For  $\hat{\beta}=(\hat{\theta},\hat{m}),$ we write ${\xi}={\theta}-\theta_0,$ $\hat{\xi}=\hat{\theta}-\theta_0,$ ${h}(Z)={m}(Z)-m_0(Z)+(\theta-\theta_0)^\top \boldsymbol{\varphi}^*(Z),$ $\hat{h}(Z)=\hat{m}(Z)-m_0(Z)+(\hat{\theta}-\theta_0)^\top \boldsymbol{\varphi}^*(Z)$ and $\tilde{X}=X-\boldsymbol{\varphi}^*(Z)$. These imply that 
\begin{equation*}
    \frac{1}{n}\sum_{i=1}^{n}\rho_{\tau}(Y_i-\theta^\top X-m(Z))=
\frac{1}{n}\sum_{i=1}^{n}\rho_{\tau}(\epsilon_i-\xi^\top \tilde{X}_i-h(Z_i)).
\end{equation*}
Denote $L_n(\xi,h)=\frac{1}{n}\sum_{i=1}^{n}\rho_{\tau}(\epsilon_i-\xi^\top \tilde{X}_i-h(Z_i)).$ We define the subgradient of the loss function $L_n$ at $\xi$ as 
\[\Psi_n(\xi,h)=\frac{\partial L_n(\xi,h)}{\partial \xi}=\mathbb{P}_n\psi_{\tau}(\xi,h),\]
where $\psi_{\tau}(\xi,h)=-\{\tau-1(\epsilon-\xi^\top \tilde{X}-h(Z)<0)\}\tilde{X}.$ Let $\Psi_0(\xi,h)=\mathbb{E}\psi_{\tau}(\xi,h).$ With $\mathcal{A}_{\delta}$ defined in \eqref{A_delta}, we further define $\tilde{\mathcal{A}}_{\delta}=\{(\xi,h)~|~\xi={\theta}-\theta_0,{h}(Z)={m}(Z)-m_0(Z)+(\theta-\theta_0)^\top \boldsymbol{\varphi}^*(Z),(\theta,m)\in\mathcal{A}_{\delta}\}$ and $\mathcal{C}_{\delta}=\{\psi_{\tau}(\xi,h)-\psi_{\tau}(\xi_0,h_0)~|~(\xi,h)\in \tilde{\mathcal{A}}_{\delta} ~\text{and}~ (\xi_0,h_0)=(0,0) \}.$
Then by analogy to the proof of Theorem \ref{Thm: convergency rate}, we have, for any $\delta>0,$ 
\begin{equation*}
    J_{[~]}(\delta,\mathcal{C}_{\delta})=\int_{0}^{\delta}\sqrt{1+ \mathcal{N}_{[~]}(\epsilon,\mathcal{C}_{\delta},L^{2}(P))}d\epsilon\lesssim \delta\sqrt{s\log\frac{R}{\delta}}.
\end{equation*}
Thus it follows
\begin{equation*}
\begin{aligned}
   &\mathbb{E}^*\{ \mathop{\sup}_{(\xi,h)\in\mathcal{C}_\delta}\big|\sqrt{n}[(\Psi_n-\Psi_0)(\xi,h)-(\Psi_n-\Psi_0)(\xi_0,h_0)]\big|\} \\=& \mathbb{E}^*\{ \mathop{\sup}_{(\xi,h)\in\mathcal{C}_\delta}\big|\sqrt{n}(\mathbb{P}_n-\mathbb{P})[\psi_{\tau}(\xi,h)-\psi_{\tau}(\xi_0,h_0)]\big|\}\\
   \lesssim& J_{[~]}(\delta,\mathcal{C}_{\delta})\Big\{ \frac{J_{[~]}(\delta,\mathcal{C}_{\delta})}{\delta^2\sqrt{n}}+1\Big\}\\
   =&\phi_n(\delta).
   \end{aligned}
\end{equation*}
This implies that 
\begin{equation*}
    \Big|\sqrt{n}[(\Psi_n-\Psi_0)(\xi,h)|_{(\xi,h)=(\hat{\xi},\hat{h})}-(\Psi_n-\Psi_0)(\xi_0,h_0)]\Big|\lesssim \phi_n(\delta)=o_p(1),
\end{equation*}
or, written alternatively,
\begin{equation}\label{proof: raw score}
\sqrt{n}\{\Psi_0(\xi,h)|_{(\xi,h)=(\hat{\xi},\hat{h})}+\Psi_n(\xi_0,h_0)\}=\sqrt{n}\{\Psi_n(\xi,h)|_{(\xi,h)=(\hat{\xi},\hat{h})}+\Psi_0(\xi_0,h_0)\} + o_p(1).
\end{equation}
Let $\tilde{Y}_i=\epsilon_i-\hat{h}(Z_i),~i=1,\ldots,n.$ Then $\hat{\xi}$ is the minimizer of $L_n^*(\xi)=\frac{1}{n}\sum_{i=1}^{n}\rho_{\tau}(\tilde{Y}_i-\xi^\top \tilde{X}_i)$ with respect to $\xi$ and
\begin{equation*}
    \Psi_n(\xi,h)\big|_{(\xi,h)=(\hat{\xi},\hat{h})}=\frac{d L_n^*(\xi)}{d \xi}\Big|_{\xi=\hat{\xi}}=\frac{1}{n}\sum_{i=1}^{n}-\{\tau-1(\tilde{Y}_i- \hat{\xi}^\top \tilde{X}_i<0)\}\tilde{X}_i.
\end{equation*}
Since $L_n^{*}$ is a continuous piecewise function of $\xi,$ it follows that the subgradient is bounded by the difference between the right and left derivatives. Thus,
\begin{equation*}
\begin{aligned}
   \Big|\frac{d L_n^*(\xi)}{d \xi}\big|_{\xi=\hat{\xi}}\Big|
   &\le\frac{2}{n}\sum_{i=1}^{n}1(\tilde{Y}_i=\hat{\xi}^\top \tilde{X}_i)|\tilde{X}_i|\\
   &\le\Big\{2\sum_{i=1}^{n}1(\tilde{Y}_i=\hat{\xi}^\top \tilde{X}_i)\Big\}\max_{i=1,\ldots,n}\Big(\frac{|\tilde{X}_i|}{n}\Big)\\
   &=o_p(\frac{1}{\sqrt{n}}),
\end{aligned}
\end{equation*}
where $|\cdot|$ and $\max(\cdot)$ operate component-wise on vector and the last equality holds due to Assumption (A\ref{assump: covariates}), (A\ref{assump: AsympNormal}) and the fact $\sum_{i=1}^{n}1(\tilde{Y}_i=\hat{\xi}^\top \tilde{X}_i)\le p$. 
Moreover, a calculation yields $\Psi_0(\xi_0,h_0)=0$, so  the left hand side of \eqref{proof: raw score} satisfies 
\begin{equation*}
\sqrt{n}\{\Psi_0(\xi,h)|_{(\xi,h)=(\hat{\xi},\hat{h})}+\Psi_n(\xi_0,h_0)\}=o_p(1),
\end{equation*}
or equivalently,
\begin{equation*}
\sqrt{n}\Psi_0(\xi,h)|_{(\xi,h)=(\hat{\xi},\hat{h})}=-\sqrt{n}\Psi_n(\xi_0,h_0)+o_p(1).
\end{equation*}

On the other hand, applying the Taylor's expansion for $\Psi_0(\xi,h)|_{(\xi,h)=(\hat{\xi},\hat{h})}$ at $(\xi_0,h_0)$, we obtain
\begin{equation*}
    \Psi_0(\xi,h)|_{(\xi,h)=(\hat{\xi},\hat{h})} = 2\mathbb{E}\{f_{\tau}(0|U)\tilde{X}\tilde{X}^\top\} (\hat{\xi}-\xi_0)+O(d^2(\hat{\beta},\beta_0)).
\end{equation*}
Here the derivative with respect to $h$ is based on the derivative of some smooth curve $\{h_{(t)}: t\in\mathbb{R},h_{(0)}=h_{0}~\text{and}~h_{(1)}=\hat{h}\}$ with respect to $t.$ Since $\hat{\xi}-\xi_0=\hat{\theta}-\theta_0$ and $\bar{\gamma}_{\bar{k}}>\bar{d}_{\bar{k}}/2$.   It follows that 
\begin{equation*}
    \sqrt{n}(\hat{\theta}-\theta_0)= \frac{1}{2}[\mathbb{E}\{f_{\tau}(0|U)\tilde{X}\tilde{X}^\top\}]^{-1}\sqrt{n}\Psi_n(\xi_0,h_0)+o_{p}(1)\rightarrow N(0,\Sigma_2^{-1}\Sigma_1\Sigma_2^{-1}).
\end{equation*}
Therefore, the result follows.
\bibliographystyle{qg}
\bibliography{DeepQR.bib} 
\end{document}